\documentclass[11pt]{amsart}
\usepackage{latexsym}
\usepackage{amsxtra}
\usepackage{amssymb}
\usepackage{amsthm}
\usepackage{xspace}

\usepackage{euscript}
\input xy
\xyoption{all} \CompileMatrices

\newtheorem{theorem}{Theorem}
\newtheorem{proposition}[theorem]{Proposition}
\newtheorem{corollary}[theorem]{Corollary}

\newtheorem{problem}[theorem]{Problem}

\begin{document}
\sloppy

\def\XX{{\bf X}}
\def\YY{{\bf Y}}
\def\EE{{\bf E}}
\def\Int{{\operatorname{Int}}}
\def\rad{{\operatorname{rad}}}
\def\Spec{\operatorname{Spec}}
\def\Sym{\operatorname{Sym}}
\def\Hom{\operatorname{Hom}}
\def\Aut{{\operatorname{Aut}}}
\def\Mor{{\operatorname{Mor}}}
\def\eval{{\operatorname{eval}}}
\def\ZZ{{\mathbb Z}}
\def\Zhat{{\widehat{\mathbb Z}}}
\def\CC{{\mathbb C}}
\def\PP{{\mathbb P}}
\def\JJ{{\mathbb J}}
\def\NN{{\mathbb N}}
\def\RR{{\mathbb R}}
\def\QQ{{\mathbb Q}}
\def\FF{{\mathbb F}}
\def\mm{{\mathfrak m}}
\def\nn{{\mathfrak n}}
\def\jj{{\mathfrak j}}
\def\aaa{{\mathfrak a}}
\def\bbb{{\mathfrak b}}
\def\ppp{{\mathfrak p}}
\def\qqq{{\mathfrak q}}
\def\PPP{{\mathfrak P}}
\def\MM{{\mathfrak M}}
\def\BB{{\mathfrak B}}
\def\jj{{\mathfrak J}}
\def\LL{{\mathfrak L}}
\def\qq{{\mathfrak Q}}
\def\rr{{\mathfrak R}}
\def\DD{{\mathfrak D}}
\def\cc{{\mathfrak S}}
\def\TT{{\mathcal{T}}}
\def\SS{{\mathcal S}}
\def\UU{{\mathcal U}}
\def\AA{{\mathcal A}}
\def\aa{{\underline{a}}}
\def\bb{{\underline{b}}}
\def\rad{\operatorname{rad}}
\def\End{\operatorname{End}}
\def\id{\operatorname{id}}
\def\mod{\operatorname{mod}}
\def\im{\operatorname{im}}
\def\ker{\operatorname{ker}}
\def\coker{\operatorname{coker}}
\def\tors{\operatorname{tors}}
\def\ord{\operatorname{ord}}
\def\Bin{\operatorname{Bin}}
\def\ww{\textup{w}}

\title{Birings and plethories of integer-valued polynomials}
\date{\today} \author{Jesse Elliott} \address{Department of Mathematics\\ California
State University, Channel Islands\\ Camarillo, California 93012}
\email{jesse.elliott@csuci.edu}

\maketitle

\begin{abstract}

Let $A$ and $B$ be commutative rings with identity.  An {\it $A$-$B$-biring} is an $A$-algebra $S$ together with a lift of the functor $\Hom_A(S,-)$ from $A$-algebras to sets to a functor from $A$-algebras to $B$-algebras.   An {\it $A$-plethory} is a monoid object in the monoidal category, equipped with the composition product, of $A$-$A$-birings.   The polynomial ring $A[X]$ is an initial object in the category of such structures.  The $D$-algebra $\Int(D)$ has such a structure if $D = A$ is a domain such that the natural $D$-algebra homomorphism $\theta_n: {\bigotimes_D}_{i = 1}^n \Int(D) \longrightarrow \Int(D^n)$ is an isomorphism for $n = 2$ and injective for $n \leq 4$.  This holds in particular if $\theta_n$ is an isomorphism for all $n$, which in turn holds, for example, if $D$ is a Krull domain or more generally a TV PVMD.  In these cases we also examine properties of the functor $\Hom_D(\Int(D),-)$ from $D$-algebras to $D$-algebras, which we hope to show is a new object worthy of investigation in the theory of integer-valued polynomials.

\ \\ 


\end{abstract}

\section{Introduction}

This paper is a summary of the results contained in the forthcoming paper \cite{ell0}.   Throughout this paper all rings and algebras are assumed commutative with identity.

For any integral domain $D$ with quotient field $K$, any set $\XX$, and any subset $\EE$ of $K^\XX$, the ring of {\it integer-valued polynomials on $\EE$ over $D$} is the subring
$$\Int(\EE,D) = \{f(\XX) \in K[\XX]: f(\EE) \subseteq D\}$$
of the polynomial ring $K[\XX]$.  In other words, $\Int(\EE,D)$ is the pullback of the direct product $D^\EE$ along the $K$-algebra homomorphism $K[\XX] \longrightarrow K^\EE$ acting by $f \longmapsto (f(\aa))_{\aa \in \EE}$.  One writes $\Int(D^\XX) = \Int(D^\XX,D)$ and $\Int(D) = \Int(D,D)$.  One also writes $\Int(D^n) = \Int(D^\XX)$ if $\XX$ is a set of cardinality $n$.

Much of the theory of integer-valued polynomial rings developed in attempts to generalize results known about $\Int(\ZZ)$ to $\Int(D)$.  This paper is concerned with finding such a generalization of a particular result about $\Int(\ZZ)$.  To state this result we need a few definitions.

A ring $A$ is said to be {\it binomial} if $A$ is $\ZZ$-torsion-free and $\frac{a(a-1)(a-2)\cdots(a-n+1)}{n!} \in A \otimes_\ZZ \QQ$ lies in $A$ for all $a \in A$ and all positive integers $n$.  For any set $\XX$ the ring $\Int(\ZZ^\XX)$ is the free binomial ring generated by $\XX$, and a $\ZZ$-torsion-free ring $A$ is binomial if and only if, for every $a \in A$, there exists a ring homomorphism $\Int(\ZZ) \longrightarrow A$ sending $X$ to $a$ \cite{ell}.  By the category of binomial rings we will mean the full subcategory of the category of rings whose objects are the binomial rings.  By \cite[Section 46]{ber} and \cite[Theorem 9.1]{ell} we have the following.

\begin{proposition}\label{introprop}
There is a functor $\Bin$ from rings to binomial rings that is left-represented by $\Int(\ZZ)$ and is a right adjoint for the inclusion from binomial rings to rings.
\end{proposition}

Our motivating problem is to generalize the above result to $\Int(D)$ for further domains $D$.   More specifically, we are interested in the following.

\begin{problem}\label{mainproblem}
Determine all domains $D$ for which $\Int(D)$ left-represents a right adjoint for the inclusion from ${\mathcal C}$ to  the category of $D$-algebras for some full subcategory ${\mathcal C}$ of the category of $D$-algebras.
\end{problem}

In particular, if $D$ is such a domain, then the functor $\Hom_D(\Int(D),-)$ from $D$-algebras to sets must lift to a functor from $D$-algebras to $D$-algebras in ${\mathcal C}$.  If $D = \ZZ$, then by Proposition \ref{introprop} this holds for the category ${\mathcal C}$ of binomial rings.   Given a domain $D$, a natural candidate for the category ${\mathcal C}$ is the category of $D$-torsion-free $D$-algebras that are ``weakly polynomially complete'' \cite[Section 7]{ell2}, where a $D$-algebra $A$ is said to be {\it weakly polynomially complete}, or {\it WPC}, if for every $a \in A$ there exists a $D$-algebra homomorphism $\Int(D) \longrightarrow A$ sending $X$ to $a$.  A binomial ring is equivalently a $\ZZ$-torsion-free WPC $\ZZ$-algebra, and for any domain $D$ the $D$-algebra $\Int(D)$ is itself WPC.  Our goal, then, is to construct a right adjoint for the inclusion from the category of $D$-torsion-free WPC $D$-algebras to the category of $D$-algebras that is left-represented by $\Int(D)$.  In our efforts to do so we found it necessary to utilize the notions of a {\it biring} and a {\it plethory}.
  
Let $A$ and $B$ be rings.  An {\it $A$-$B$-biring} is an $A$-algebra $S$ together with the structure on $S$ of a $B$-algebra object in the opposite category of the category of $A$-algebras.  Thus an $A$-$B$-biring is an $A$-algebra $S$ equipped with two {\it binary co-operations} $S \longmapsto S \otimes_A S$, called {\it co-addition} and {\it co-multiplication} (both of which are $A$-algebra homomorphisms), along with a {\it co-$B$-linear structure} $B \longrightarrow \Hom_A(S,A)$, satisfying laws dual to those defining the $A$-algebras.  By Yoneda's lemma, an $A$-$B$-biring is equivalently an $A$-algebra $S$ together with a lift of the covariant functor $\Hom_A(S,-)$ it represents to a functor from the category of $A$-algebras to the category of $B$-algebras.  (See any of \cite{ber,bor,tall} for the details.)  For example, the polynomial ring $A[X]$ is an $A$-$A$-biring as it represents the identity functor from the category of $A$-algebras to itself.  Co-addition acts by $X \longmapsto X \otimes 1 + 1 \otimes X$, co-multiplication by $X \longmapsto X \otimes X$, and the co-linear structure by $a \longmapsto (f \longmapsto f(a))$.

\begin{proposition}\label{biringprop}
Let $D$ be an integral domain.
\begin{enumerate}
\item The existence of a $D$-$D$-biring structure on $\Int(D)$ is equivalent to the existence of a lift of the functor $\Hom_D(\Int(D),-)$ from $D$-algebras to sets to a functor from $D$-algebras to $D$-algebras.
\item A $D$-$D$-biring structure on $\Int(D)$ is compatible with the $D$-$D$-biring structure on $D[X]$, that is, the inclusion $D[X] \longrightarrow \Int(D)$ is a homomorphism of $D$-$D$-birings, if and only if the natural map $\Hom_D(\Int(D),A) \longrightarrow A$ given by $\varphi \longmapsto \varphi(X)$ is a $D$-algebra homomorphism for every $D$-algebra $A$.
\end{enumerate}
\end{proposition}

Consequently, any solution to Problem \ref{mainproblem} would yield conditions on integral domains $D$ under which the $D$-algebra $\Int(D)$ has a $D$-$D$-biring structure.  Regarding the latter problem we have the following.

\begin{theorem}\label{th} Assume that $\Int(D)$ is flat over $D$, or more generally that the $n$-th tensor power $\Int(D)^{\otimes n}$ of $\Int(D)$ over $D$ is $D$-torsion-free for $n \leq 4$.  Then the domain $\Int(D)$ has a (necessarily unique) $D$-$D$-biring structure that is compatible with the $D$-$D$-biring structure on $D[X]$ if and only if for every $f \in \Int(D)$ the polynomials $f(X+Y)$ and $f(XY)$ both can be expressed as sums of polynomials of the form $g(X)h(Y)$ for $g,h \in \Int(D)$.
\end{theorem}  

In analogy with ordinary polynomial rings, there is for any set $\XX$ a canonical $D$-algebra homomorphism $$\theta_\XX: \bigotimes_{X \in \XX} \Int(D) \longrightarrow \Int(D^\XX),$$ where the (possibly infinite) tensor product is over $D$ and is a coproduct in the category of $D$-algebras.  However, we do not know whether or not $\theta_\XX$ is an isomorphism for every domain $D$ and every set $\XX$.  There are several large classes of domains for which $\theta_\XX$ is an isomorphism for all $\XX$, such as the Krull domains, the almost Newtonian domains \cite[Section 5]{ell2}, and the PVMDs $D$ such that $\Int(D_\mm) = \Int(D)_\mm$ for every maximal ideal $\mm$, hence the TV PVMDs as well.  (See \cite{hwa} for the definition of a PVMD and a TV PVMD.)  We say that the domain $D$ is {\it polynomially composite} if $\theta_\XX$ is an isomorphism for every set $\XX$.  Section 4 of \cite{ell3} collects several known classes of polynomial composite domains.  

\begin{corollary}\label{polycomp}
If $D$ is a polynomially composite domain, and in particular if $D$ is a Krull domain or TV PVMD, then $\Int(D)$ has a unique $D$-$D$-biring structure such that the inclusion $D[X] \longrightarrow \Int(D)$ is a homomorphism of $D$-$D$-birings.
\end{corollary}

By \cite[Proposition 9.3]{ell} one has $\Bin(A) \cong \ZZ_p$ for any integral domain $A$ of characteristic $p$, where $\ZZ_p$ denotes the ring of $p$-adic integers, and in particular one has $\Bin(\FF_p) \cong \ZZ_p$. This generalizes as follows.

\begin{proposition}
Let $D$ be a Dedekind domain, and let $\mm$ be a maximal ideal of $D$ with finite residue field.  Then the map
$$\widehat{D}_\mm \longrightarrow \Hom_D(\Int(D), D/\mm)$$
acting by $\alpha \longmapsto (f \longmapsto f(\alpha) \mod \mm \widehat{D}_\mm)$
is a $D$-algebra isomorphism.  More generally, for any domain extension $A$ of $D$ with $\mm A = 0$, the diagram
\begin{eqnarray*}
\SelectTips{cm}{11}\xymatrix{
 {\widehat{D}_\mm} \ar[r] \ar[dr] & {\Hom_D(\Int(D), D/\mm)} \ar[d] \\
 & {\Hom_D(\Int(D), A)}}
\end{eqnarray*}
is a commutative diagram of $D$-algebra isomorphisms.
\end{proposition}

By \cite[Proposition 1.4]{bor}, for any $A$-$B$-biring $S$, the lifted functor $\Hom_A(S,-)$ from $A$-algebras to $B$-algebras has a left adjoint, denoted $S \odot_A -$.  In analogy with the tensor product,  the $A$-algebra $S \odot_A R$ for any $B$-algebra $R$ is the $A$-algebra generated by the symbols $s \odot r$ for all $s \in S$ and $r \in R$, subject to the relations \cite[1.3.1--2]{bor}.  If $S$ and $T$ are $A$-$A$-birings, then so is $S \odot_A T$, and the category of $A$-$A$-birings equipped with the operation $\odot_A$ is monoidal with unit $A[X]$.   An {\it $A$-plethory} is a monoid object in that monoidal category, that is, it is an $A$-$A$-biring $P$ together with an associative map $\circ : P \odot_A P \longrightarrow P$ of $A$-$A$-birings (called {\it composition}) possessing a unit $e: A[X] \longrightarrow P$.  (See any of \cite{ber,bor,tall} for details on these constructions.)  An $A$-plethory is also known as an {\it $A$-$A$-biring monad object}, an {\it $A$-$A$-biring triple}, or a {\it Tall-Wraith monad object} in the category of $A$-algebras.  For example, for any ring $A$, the polynomial ring $A[X]$ has the structure of an $A$-plethory and in fact is an initial object in the category of $A$-plethories.

\begin{proposition}
Let $D$ be an integral domain.  Any $D$-$D$-biring structure on $\Int(D)$ compatible with that on $D[X]$ extends uniquely to a $D$-plethory structure on $\Int(D)$ with unit given by the inclusion $D[X] \longrightarrow \Int(D)$.  Composition $\circ: \Int(D) \odot_D \Int(D) \longrightarrow \Int(D)$ acts by ordinary composition on elements of the form $f \odot g$, that is, one has $\circ: f \odot g \longmapsto f \circ g$ for all $f, g \in \Int(D)$.
\end{proposition}

The following theorem gives a partial solution to Problem \ref{mainproblem}.

\begin{theorem}\label{introthm}
Let $D$ be an integral domain.
\begin{enumerate}
\item Assume that there exists a $D$-$D$-biring structure on $\Int(D)$ compatible with that on $D[X]$.  Then the functors $\Hom_D(\Int(D),-)$ and $\Int(D) \odot_D -$ are right and left adjoints, respectively, of the inclusion from $D$-torsion-free WPC algebras to $D$-algebras if and only if the $D$-algebras $\Hom_D(\Int(D),A)$ and $\Int(D) \odot_D A$ are $D$-torsion-free for any $D$-algebra $A$.
\item If $D$ is a PID with all residue fields finite, then the hypotheses (and therefore the conclusion) of statement (1) hold.
\end{enumerate}
\end{theorem}

Note that if $\Int(D) = D[X]$, which, for example, holds by \cite[Corollary I.3.7]{cah} if $D$ has no finite residue fields, then $\Hom_D(\Int(D), A)$ is naturally isomorphic to $A$ and in particular is not $D$-torsion-free if $A$ is not $D$-torsion-free.  Of course in that case every $D$-algebra is WPC.

\section{WPC $D$-algebras and tensor powers of $\Int(D)$}

As in \cite[Section 7]{ell2} and as in the introduction, we will say that a $D$-algebra $A$ is {\it weakly polynomially complete}, or {\it WPC}, if for every $a \in A$ there exists a $D$-algebra homomorphism $\Int(D) \longrightarrow A$ sending $X$ to $a$. 
A $D$-torsion-free $D$-algebra $A$ is WPC if and only if $f(A) \subseteq A$ for all $f \in \Int(D) \subseteq (A \otimes_D K)[X]$.  In particular, a domain extension $A$ of $D$ is WPC if and only if $\Int(D) \subseteq \Int(A)$. 

For any set $\XX$, the smallest subring of $\Int(D^\XX)$ containing $D[\XX]$ that is closed under precomposition by elements of $\Int(D)$ is denoted $\Int_\ww(D^\XX)$.  For any domain $D$ (finite or infinite), the domain $\Int_\ww(D^\XX)$ is the free WPC extension of $D$ generated by $\XX$ \cite[Proposition 7.2]{ell2}.  It is also the {\it weak polynomial completion $w_D(D[\XX])$ of $D[\XX]$ with respect to $D$}, as defined in \cite[Section 8]{ell2} and in Proposition \ref{weakcompositum} below.

If $\Int_\otimes(D^\XX)$ denotes the image of the $D$-algebra homomorphism $\theta_\XX: \bigotimes_{X \in \XX} \Int(D) \longrightarrow \Int(D^\XX)$, then we have $\Int_\otimes(D^\XX) \subseteq \Int_\ww(D^\XX)$, and equality holds for a given set $\XX$ if and only if $\Int_\otimes(D^\XX)$ is a WPC extension of $D$.  If equality holds for any set $\XX$ then we will say that $D$ is {\it weakly polynomially composite}.

\begin{proposition}\label{manyvarprop2}  The following conditions are equivalent for any integral domain $D$.
\begin{enumerate}
\item $D$ is weakly polynomially composite.
\item $\Int_\otimes(D^\XX)$ is a WPC extension of $D$ for any set $\XX$.
\item $\Int_\otimes(D^n)$ is a WPC extension of $D$ for some integer $n > 1$.
\item $\Int_\otimes(D^2)$ is a WPC extension of $D$.
\item For any element $f$ of $\Int(D)$, the polynomials $f(X+Y)$ and $f(XY)$ lie in the image of $\theta_{\{X,Y\}}$.
\item The compositum of any collection of  WPC $D$-algebras of $D$ contained in some $D$-torsion-free $D$-algebra is again a WPC $D$-algebra.
\item The compositum of any collection of WPC extensions of $D$ contained in some domain extension of $D$ is again a WPC extension of $D$.
\end{enumerate}
\end{proposition}

Clearly polynomial compositeness implies weak polynomial compositeness. 

At the end of Section 8 of \cite{ell2} it is noted how to construct the left adjoint of the inclusion functor from WPC domain extensions of $D$ to domain extensions of $D$.  The proof can be easily generalized to establish the following.

\begin{proposition}\label{weakcompletion}
Let $D$ be a domain with quotient field $K$, and let $A$ be a $D$-torsion-free $D$-algebra.
\begin{enumerate}
\item $A$ is contained in a smallest $D$-torsion-free WPC $D$-algebra, denoted $w_D(A)$, equal to the intersection of all WPC $D$-algebras containing $A$ and contained in $A \otimes_D K$.
\item One has $w_D(A) = A$ if and only if $A$ is WPC, and $w_D(A)$ is a domain if and only if $A$ is a domain.
\item One has $w_D(A) \cong \Int_\ww(D^\XX)/((\ker \varphi)K \cap \Int_\ww(D^\XX))$ for any surjective $D$-algebra homomorphism $\varphi: D[\XX] \longrightarrow A$.
\item The association $A \longmapsto w_D(A)$ defines a functor from the category of $D$-torsion-free $D$-algebras to the category of $D$-torsion-free WPC $D$-algebras---both categories with morphisms as $D$-algebra homomorphisms---that is a left adjoint for the inclusion functor.
\end{enumerate}
\end{proposition}

Assuming that $D$ is weakly polynomially composite, we can also construct the right adjoint of the inclusion functor from $D$-torsion-free WPC $D$-algebras to $D$-torsion-free $D$-algebras.

\begin{proposition}\label{weakcompositum}
Let $D$ be a weakly polynomially composite domain, and let $A$ be a $D$-torsion-free $D$-algebra.
\begin{enumerate}
\item $A$ contains a largest WPC $D$-algebra, denoted $w^D(A)$, equal to the compositum of all WPC $D$-algebras contained in $A$.
\item One has $w^D(A) = A$ if and only if $A$ is WPC. 
\item One has $w^D(A) = \{a \in A: a = \varphi(X) \mbox{ for some } \varphi \in \Hom_D(\Int(D),A)\}$.
\item The association $A \longmapsto w^D(A)$ defines a functor from the category of $D$-torsion-free $D$-algebras to the category of $D$-torsion-free WPC $D$-algebras---both categories with morphisms as $D$-algebra homomorphisms---that is a right adjoint for the inclusion functor.
\end{enumerate}
\end{proposition}

\section{Biring and plethory structure on $\Int(D)$}

The following result implies Theorem \ref{th} and Corollary \ref{polycomp} of the introduction.

\begin{theorem}\label{biringstructure}  Let $D$ be an integral domain.
\begin{enumerate}
\item If the domain $\Int(D)$ has a $D$-$D$-biring structure such that the inclusion $D[X] \longrightarrow \Int(D)$ is a homomorphism of $D$-$D$-birings, then $D$ is weakly polynomially composite.
\item Assume that the $n$-th tensor power $\Int(D)^{\otimes n}$ of $\Int(D)$ over $D$ is $D$-torsion-free for $n \leq 4$.  Then $\Int(D)$ has a unique $D$-$D$-biring structure such that the inclusion $D[X] \longrightarrow \Int(D)$ is a homomorphism of $D$-$D$-birings if $D$ is weakly polynomially composite.
\end{enumerate}
\end{theorem}

The plethory $A[X]$ is an initial object in the category of $A$-plethories.  Like $A[X]$, and in particular like the domain $D[X]$, the domain $\Int(D)$ has its own ``internal'' operation of composition.  This leads to the following result.

\begin{proposition}\label{intplethory}
Let $D$ be an integral domain.  Any $D$-$D$-biring structure on $\Int(D)$ such that the inclusion $D[X] \longrightarrow \Int(D)$ is a homomorphism of $D$-$D$-birings extends uniquely to a $D$-plethory structure on $\Int(D)$ with unit given by the inclusion $D[X] \longrightarrow \Int(D)$.  Composition $\circ: \Int(D) \odot_D \Int(D) \longrightarrow \Int(D)$ acts by ordinary composition on elements of the form $f \odot g$, that is, one has $\circ: f \odot g \longmapsto f \circ g$ for all $f, g \in \Int(D)$.
\end{proposition}

\begin{corollary}\label{maincor}
If $D$ is a polynomially composite domain, and in particular if $D$ is a Krull domain or TV PVMD, then $\Int(D)$ has a unique $D$-plethory structure with unit given by the inclusion $D[X] \longrightarrow \Int(D)$.
\end{corollary}

Let $A$ be a ring and $P$ an $A$-plethory.  A {\it $P$-ring} is an $A$-algebra $R$ together with an $A$-algebra homomorphism $\circ : P \odot_A R \longrightarrow R$ such that $(\alpha \circ \beta) \circ r = \alpha \circ (\beta \circ r)$ and $e \circ r = e$ for all $\alpha, \beta \in P$ and all $r \in R$, where $e$ is the image of $X$ in the unit $A[X] \longrightarrow P$ \cite[1.9]{bor}.  Such a map $\circ$ is said to be a {\it left action of $P$ on $R$}.  For example, $P$ itself has a structure of a $P$-ring, as do the $A$-algebras $P \odot_A R$ and $\Hom_A(P,R)$ for any $A$-algebra $R$ \cite[1.10]{bor}, with left actions given by
\begin{eqnarray*}
P \odot_A (P \odot_A R) & \longrightarrow & P \odot_A R \\
\alpha \odot (\beta \odot r) & \longmapsto & (\alpha \circ \beta) \odot r
\end{eqnarray*}
and
\begin{eqnarray*}
P \odot_A \Hom_A(P,R) & \longrightarrow & \Hom_A(P,R) \\
\alpha \odot \varphi & \longmapsto & (\beta \longmapsto \varphi(\beta \circ \alpha)),
\end{eqnarray*}
respectively.  Moreover, the functors $P \odot_A -$ and $W_P = \Hom_A(P,-)$ from $A$-algebras to $P$-rings are left and right adjoints, respectively, for the forgetful functor from $P$-rings to $A$-algebras \cite[1.10]{bor}.

For any $A$-plethory $P$, the $P$-ring $W_P(R) = \Hom_A(P,R)$ of any $A$-algebra $R$ is called the {\it $P$-Witt ring} of $R$.  This terminology comes from the fact that, if $P$ is the $\ZZ$-plethory $\Lambda$ of \cite[2.11]{bor}, then a $P$-ring is equivalently a $\lambda$-ring, and the functor $W_P$ is isomorphic to the universal $\lambda$-ring functor $\Lambda$.   If $P$ is the $\ZZ$-plethory $\Int(\ZZ)$, then a $P$-ring is equivalently a binomial ring, and the functor $W_P$ is isomorphic to the functor $\Bin$.  The latter fact generalizes to the following result, which implies Theorem \ref{introthm} of the introduction.

\begin{theorem}\label{maintheorem}  Let $D$ be an integral domain such that $\Int(D)$ has a $D$-plethory structure with unit given by the inclusion $D[X] \longrightarrow \Int(D)$, and let $A$ be a $D$-algebra.
\begin{enumerate}
\item If there exists an $\Int(D)$-ring structure on $A$, then $A$ is WPC.
\item If $A$ is $D$-torsion-free,  then there exists a (necessarily unique) $\Int(D)$-ring structure on $A$ if and only if $A$ is WPC.
\item If $A$ is $D$-torsion-free, then the $D$-algebra homomorphism $W_{\Int(D)}(A) \longrightarrow A$ is an inclusion with image equal to $w^D(A)$, and the functor $w^D$ is therefore isomorphic to the functor $W_{\Int(D)}$ restricted to the category of $D$-torsion-free $D$-algebras.
\item If $A$ is $D$-torsion-free, then the $D$-algebra homomorphism $\Int(D) \odot_D A \longrightarrow A \otimes_D K$ acting by $f \odot a \longmapsto f(a)$ has image equal to $w_D(A)$, and the functor $w_D$ is therefore isomorphic to the functor $T \circ (\Int(D) \odot_D -)$ restricted to the category of $D$-torsion-free $D$-algebras, where $T(B)$ for any $D$-algebra $B$ denotes the image of $B$ in $B \otimes_D K$, where $K$ is the quotient field of $D$.
\item If $A$ is $D$-torsion-free and WPC, then the natural $D$-algebra homomorphisms $W_{\Int(D)}(A) \longrightarrow A$ and $A \longrightarrow T(\Int(D) \odot_D A)$ are isomorphisms.
\item The functor $T \circ (\Int(D) \odot_D -)$ is a left adjoint for the inclusion from $D$-torsion-free WPC $D$-algebras to $D$-algebras.
\item The functors $W_{\Int(D)}$ and $\Int(D) \odot_D-$ are right and left adjoints, respectively, for the inclusion from $D$-torsion-free WPC $D$-algebras to $D$-algebras if and only if the $\Int(D)$-rings $W_{\Int(D)}(R)$ and $\Int(D) \odot_D R$ are $D$-torsion-free for every $D$-algebra $R$.
\item Every $\Int(D)$-ring is $D$-torsion-free if $D$ is a PID with finite residue fields.
\end{enumerate}
\end{theorem}

We end with the following problem.

\begin{problem}
Determine equivalent conditions on an integral domain $D$ so that the $D$-algebra $\Int(D)$ has a $D$-plethory structure with unit given by the inclusion $D[X] \longrightarrow \Int(D)$ and so that the $D$-algebras $W_{\Int(D)}(R)$ and $\Int(D) \odot_D R$ are $D$-torsion-free for every $D$-algebra $R$.
\end{problem}


\begin{thebibliography}{}

\bibitem{ber} G.\ Bergman and A. Hausknecht, {\it Cogroups and Co-Rings in Categories of Associative Rings}, Mathematical Surveys and Monographs, Volume 45, American Mathematical Society, 1996.

\bibitem{bor} J.\ Borger and B.\ Wieland, Plethystic algebra, Adv. Math. 194 (2005) 246--283.

\bibitem{cah} P.-J.\ Cahen and J.-L.\ Chabert, {\it Integer-Valued Polynomials},  Mathematical Surveys and Monographs, vol.\ 48, American Mathematical Society, 1997.

\bibitem{ell} J. Elliott, Binomial rings, integer-valued polynomials, and $\lambda$-rings,
J.\ Pure Appl.\ Alg.\ 207 (2006) 165--185.

\bibitem{ell2} J.\ Elliott, Universal properties of integer-valued polynomial rings, J.\ Algebra \ 318 (2007) 68--92.

\bibitem{ell3} J.\ Elliott, Some new approaches to integer-valued polynomial rings, in {\it Commutative Algebra and its Applications: Proceedings of the Fifth Interational Fez Conference on Commutative Algebra and Applications}, Eds.\ Fontana, Kabbaj, Olberding, and Swanson, de Gruyter, New York, 2009.

\bibitem{ell0} J.\ Elliott, Biring and plethory structures on integer-valued polynomial rings, to be submitted for publication.

\bibitem{hwa} C.\ J.\ Hwang and G.\ W.\ Chang, Bull.\ Korean Math.\ Soc.\ 35 (2) (1998) 259--268.

\bibitem{tall} D.\ O.\ Tall and G.\ C.\ Wraith, Representable functors and operations on rings, Proc. London Math Soc. (3) 20 (1970) 619--643.

\end{thebibliography}
\end{document}